\documentclass[12pt]{amsart}
\usepackage{amsmath,amssymb,amsbsy,amsfonts,latexsym,amsopn,amstext,
                                               amsxtra,euscript,amscd,bm,cite}
\usepackage{url}
\usepackage[colorlinks,linkcolor=blue,anchorcolor=blue,citecolor=blue,backref=page]{hyperref}
\usepackage{color}
\usepackage{float}

\usepackage{comment}
\begin{document}

\newtheorem{theorem}{Theorem}
\newtheorem{lemma}[theorem]{Lemma}
\newtheorem{claim}[theorem]{Claim}
\newtheorem{cor}[theorem]{Corollary}
\newtheorem{prop}[theorem]{Proposition}
\newtheorem{definition}{Definition}
\newtheorem{question}[theorem]{Question}
\newtheorem{remark}[theorem]{Remark}
\newcommand{\hh}{{{\mathrm h}}}

\numberwithin{equation}{section}
\numberwithin{theorem}{section}
% \numberwithin{remark}{section}
\numberwithin{table}{section}

\title[Large sieve inequality for power moduli]{Large sieve inequality for power moduli}
\date{\today}

\author[M. Munsch]{Marc Munsch}
\address{5010 Institut f\"{u}r Analysis und Zahlentheorie
8010 Graz, Steyrergasse 30, Graz}
\email{munsch@math.tugraz.at}

\footnotetext{
2010 Mathematics Subject Classification. 
Primary  Secondary : }

\begin{abstract}
In this note we give a new bound for large sieve with characters to power moduli which improves in some range of the parameters the previous bounds of Baier/Zhao and Halupczok.
\end{abstract}

\bibliographystyle{plain}
\maketitle

\section{Introduction}

 The large sieve originated in the work of Linnik in the early forties and became a fundamental tool in number theory. Since then, it has been a topic largely studied and many variants were obtained. The case of power moduli received a particular interest due to its applications. Let $\{a_n\}$ denote a sequence of complex numbers, $M,N,k$ positive integers and let $Q$ be a real number $\geq 1$. The main goal is to obtain an estimate of the following kind

\begin{equation}\label{largesieve} \sum_{q \leq Q}\sum_{a=1 \atop (a,q)=1}^{q^k} \left\vert \sum_{n=M+1}^{M+N} a_n e\left(\frac{a}{q^k}n\right)\right\vert^2 \ll \Delta \sum_{n=M+1}^{M+N} \vert a_n\vert^2\end{equation} where $e(\alpha):=\exp(2i\pi \alpha)$ for $\alpha \in \mathbb{R}$. Let us recall the general large sieve inequality. We say that a set of real numbers $\{x_k\}$ is $\delta$-spaced modulo $1$ if $\|x_k-x_j\| \geq \delta$ for all $k\neq j$  where $\|x\|$ denotes the distance of a real number $x$ to its closest integer.

\begin{theorem}\cite[Theorem $2$, Chapter $27$]{Dav} Let $\{a_n\}$ be a sequence of complex numbers, $\{x_k\}$ be a set of real numbers which is $\delta$-spaced modulo $1$, and $M,N$ be integers. Then we have 

\begin{equation}\label{classic} \sum_{k}\left\vert \sum_{n=M+1}^{M+N} a_n e(x_k n)\right\vert^2 \ll (\delta^{-1}+N)  \sum_{n=M+1}^{M+N} \vert a_n\vert^2. \end{equation}
 
\end{theorem}As it is for instance pointed out by Zhao in \cite{Zhaoacta}, the classical large sieve inequality \eqref{classic} enables us to obtain \eqref{largesieve} with \begin{equation}\label{trivialdelta} \Delta= \min\left\{Q^{2k}+N,Q(Q^k+N)\right\}.\end{equation} Furthermore we can expect the fractions with power denominator to be regularly spaced. This led Zhao to conjecture that we can take $\Delta=Q^{\varepsilon}(Q^{k+1}+N)$ in \eqref{largesieve}. The bounds coming from \eqref{trivialdelta}  imply the conjecture of Zhao in the ranges $Q \leq N^{1/2k}$ and $Q \geq N^{1/k}$. \\ %In a recent work, Kerr \cite{KerrB} investigates Zhao's conjecture and obtains a $l^1 \rightarrow l^2$ large sieve inequality. 

Several authors obtained improvements of \eqref{trivialdelta} in the critical range $Q^k \leq N \leq Q^{2k}$. First Zhao proved an inequality of type \eqref{largesieve} with 

\begin{equation}\label{Zhaobound} \Delta = \left[Q^{k+1}+\left(NQ^{\frac{\kappa-1}{\kappa}}+N^{1-\kappa}Q^{\frac{\kappa+k}{\kappa}}\right)N^{\varepsilon}\right] \end{equation} where $\kappa=2^{k-1}$. Baier and Zhao proved in \cite{BZhaoIJNT} that we can take 
\begin{equation}\label{BaierZhao} \Delta= (Q^{k+1}+N+N^{1/2+\varepsilon}Q^k) (\log\log 10 NQ)^{k+1}\end{equation} which improves  \eqref{Zhaobound} in the range $N^{\frac{1}{2k}+\varepsilon} \ll Q \ll N^{\frac{(\kappa-2)}{(2(\kappa-1)\kappa-2k)}-\varepsilon}$. In particular, for $k=3$, this led to an improvement in the range $N^{1/6+\varepsilon} \ll Q \ll N^{1/5-\varepsilon}$. Using a Fourier analytic method, they obtained a further improvement for $k=3$. These results have been sharpened by Halupczok \cite{KarinIJNT} using the breakthrough work of Wooley on Vinogradov mean value conjecture \cite{Wooley}. Precisely, she proved (and further generalized to any polynomial moduli \cite{Karinquart}) that we can take 
\begin{equation}\label{Karinbound} \Delta= (QN)^{\epsilon} \left(Q^{k+1}+Q^{1-\delta}N+Q^{1+k\delta}N^{1-\delta}\right)\end{equation} where $\delta=1/(2k(k-1))$. It should be noticed that this improved \eqref{Zhaobound} for $k$ sufficiently large and \eqref{BaierZhao} for all $k\geq 3$ and $Q^k \leq N\leq Q^{2k-2+2\delta}$. In the meantime, the conjecture in
Vinogradov mean value theorem was proved for all degrees exceeding 3 by Bourgain, Demeter and Guth in \cite{Bourgainvino} and later using  another method by Wooley \cite{Wooley2}. By this impressive work, the best possible exponent in Vinogradov mean's value theorem has been confirmed. Thus we can replace $\delta$ in \eqref{Karinbound} by $2\delta$, with only few obvious changes to be made in her proof (see the discussion in \cite{Karinsurvey}). To conclude \eqref{largesieve} holds with $\Delta= (QN)^{\varepsilon} A_k(Q,N)$ where 

\begin{equation}\label{KarinMvt}A_k(Q,N)=\left(Q^{k+1}+Q^{1-\frac{1}{k(k-1)}}N+Q^{1+\frac{1}{k-1}}N^{1-\frac{1}{k(k-1)}}\right) .\end{equation} \\  She has further refined the bound \eqref{KarinMvt} in \cite{preKarin} and obtained 

\begin{equation}\label{Karin2k}\Delta= Q^{\epsilon}\left(Q^{k+1},\min\left\{A_k(Q, N), N^{1-\omega} Q^{
1+(2k-1)\omega}\right\}\right)  \end{equation} with $\omega=1/((k-1)(k-2)+2)$. \\

 We propose an improvement of the existing results in some ranges of the parameters $N$ and $Q$. In order to do so, we employ an elementary method which makes use of a bound on the number of points modulo an integer of polynomial equations due to Cilleruelo, Garaev, Ostafe and Shparlinski \cite{JavierIgorboxes}. Our result is the following
 
 \begin{theorem}\label{shortrange} With $\{a_n\}$, $Q,M$, and $N$ as before such that $N^{1/2k} \leq Q \leq N^{\frac{1}{k}}$, we have 
$$\sum_{q \leq Q}\sum_{a=1 \atop (a,q)=1}^{q^k} \left\vert \sum_{n=M+1}^{M+N} a_n e\left(\frac{a}{q^k}n\right)\right\vert^2 \ll (QN)^{\varepsilon}Q^{1+1/(k+1)}N^{1-\frac{1}{k(k+1)}} \sum_{n=M+1}^{M+N} \vert a_n\vert^2.$$
   \end{theorem}

\begin{remark}This improves \eqref{BaierZhao} for all $k\geq 3$ when $Q^k \leq N\ll Q^{2k-2+\frac{2(k-2)}{k^2+k-2}}$ and improves \eqref{Karin2k} when $ Q^{k+1+\frac{2}{k-1}}\leq N \leq Q^{2k-1+O(1/k^3)}$. The range where it improves all the previous bounds becomes non empty as soon as $k\geq 4$ and covers almost the whole range except the corners. Additionally it is a step towards Conjecture $21$ raised in \cite{preKarin}. \end{remark}   

\begin{remark} This result can be easily generalized to any polynomial $f(q)$ of degree $k$ in a similar way as in \cite[Corollary $2.3$]{Karinquart}. For the sake of simplicity and coherence, we restricted the presentation to the case of monomials. \end{remark}

\section{Proof of the main result}

 To begin with, we follow a standard path and define a subset of Farey fractions
$$ \mathcal{S}(Q):=\left\{ \frac{a}{q^k}, \,\, (a,q)=1,\, 1\leq a < q^k,\, Q \leq q \leq 2Q\right\}.$$ It is easy to remark that two distinct elements of $\mathcal{S}(Q)$ are $1/Q^{2k}$ spaced. In the case of squares, Zhao \cite{Zhaoacta} introduced a quantity to measure the spacings between these Farey fractions. Similarly we define

\begin{equation}\label{numberclose} M(N,Q)= \max_{x \in \mathcal{S}(Q)}\# \left\{y \in \mathcal{S}(Q): \|x-y\| <1/2N \right\}.\end{equation} As noticed in \cite{Zhaoacta}, any good estimate on this quantity leads to an equality of type \eqref{largesieve}. We prove the following bound.

\begin{lemma}\label{estimateclose} For any $\varepsilon>0$ and integer $N$ and $Q^k \leq N \leq Q^{2k}$, we have
\begin{equation}\label{estimatefrac} M(N,Q) \ll (QN)^{\varepsilon} \left(Q^{1+1/(k+1)}N^{-\frac{1}{k(k+1)}}\right)\end{equation}
 \end{lemma} In order to do so, we use a result bounding  the number of solutions of polynomial equations in boxes. This can be proved in a similar way as Theorem $1$ of \cite{JavierIgorboxes} or could be deduced from the generalization of Kerr \cite[Theorem $3.1$]{Kerrboxes}.
\begin{theorem}\label{boxes} Let $f$ a polynomial of degree $k\geq 2$ with leading coefficient coprime to $m$, $1 \leq  H,R \leq m$ and integers $K,L$. We define by $N(H,R;K,L)$ the number of solutions to the congruence 
\begin{equation}\label{congruence} f(x) \equiv y \,\,(\bmod\, m) \end{equation} with 
\begin{equation}\label{constraint} (x,y) \in [K+1,K+H] \times [L+1,L+R]. \end{equation}
 Then, uniformly over arbitrary integers $K$ and $L$, we have for any $\varepsilon > 0$
\begin{equation}\label{JaviIgor} N(H,R;K,L) \ll H\left((R/m)^{1/j(k)+\varepsilon} + (R/H^k)^{1/2j(k)+\varepsilon}\right).  \end{equation}
 \end{theorem} At the time these articles \cite{JavierIgorboxes,Kerrboxes} were written, the authors pointed out that $j(k)=k(k+1)$ was an admissible value following the work of Wooley \cite{Wooley}. As already noticed above, the resolution of the Vinogradov mean value conjecture allows nowadays to take $j(k)=\frac{k(k+1)}{2}$. \\

\subsection{Proof of Lemma \ref{estimateclose}}Let $x=\frac{a}{q^k}$ with $(a,q)=1$ and $y=\frac{b}{r^k}$ with $(a,q)=(b,r)=1$. We want to estimate the number of pairs $(b,r)$ with $(b,r)=1$ such that
\begin{equation}\label{maj}\left\|\frac{a}{q^k}-\frac{b}{r^k}\right\|=\frac{\vert ar^k-bq^k\vert}{q^kr^k} <1/2N.  \end{equation} Setting $z=ar^k-bq^k$, our problem boils down to estimate the number of pairs $(b,r)$ such that $\vert z\vert \ll Q^{2k}/N$. Equivalently, we will count the number of pairs $(r,z)$. The number of solutions is bounded above by the number of pairs with $r\sim Q$ and $\vert z\vert \ll Q^{2k}/N$ which are solutions of the congruence 

\begin{equation}\label{cong}ar^k= z \,(\,\bmod\, q^k).\end{equation} Applying Theorem \ref{boxes} to the polynomial $f(x)=ax^k$ with parameters $H=Q$, $R=Q^{2k}/N$, and $m=q^k$ we deduce that the numbers of pairs $(r,z)$ verifying \eqref{cong} is bounded above by 

$$(QN)^{\varepsilon}Q(Q^k/N)^{\frac{1}{k(k+1)}}. $$ The result follows.

\subsection{Proof of Theorem \ref{shortrange}}

  We first split the range of $q$ into $\log Q$ dyadic intervals. Then we divide the interval $(0,1)$ into $N$ subintervals $I_k$ of size $1/N$. In each of the intervals $I_{2k}$ we pick one element of $\mathcal{S}(Q)$ (we do the same for the odd indices). To avoid any problems at the edges we split the odd indices into two subsets: the first including elements of the initial interval and the second the fractions of the last interval. This gives a sequence of elements of $\mathcal{S}(Q)$ which are $1/N$-spaced. We can therefore apply  \eqref{classic} with $\delta=1/N$. We repeat the procedure again. In order to pick all the possible fractions from $\mathcal{S}(Q)$, we have to repeat the process at most $M(N,Q)$ times and obtain $$ \sum_{q \leq Q}\sum_{a=1 \atop (a,q)=1}^{q^k} \left\vert \sum_{n=M+1}^{M+N} a_n e\left(\frac{a}{q^k}n\right)\right\vert^2 \ll M(N,Q) N \sum_{n=M+1}^{M+N} \vert a_n\vert^2.$$ Using Lemma \ref{estimateclose}, we conclude the proof.

 %Using the duality principle and Poisson summation formula (as it is done in the proof of \cite[Theorem $2$]{Zhaoacta}), any estimate on \eqref{numberclose} of size $M$ directly transfers to an estimate  $\Delta \ll MN$ in \eqref{largesieve}. This concludes the proof of Theorem \ref{shortrange}.

%\textbf{Question:} Estimate the following quantity
%$$\sup_{a/q \in \mathcal{F}(Q)}\left\vert \left\{(a_1,q_1), \vert a/q^2-a_1/q_1^2 \vert \leq 1/N^3 \right\}\right\vert$$
%Interesting case: $N=Q^{1/3}$ applications to large sieve Zhao conjecture, then several classic applications plus if good enough application to Balog-Rusza bound on exponential sums over squarefree numbers, other interesting potential consequences? \\

% Ideas: Pell equations  \\
 
 % The problem boils down to count the following: for any couple $(a,q)$ with $1\leq a \leq q^2$ and $Q \leq q \leq 2Q$:
  
 % $$ Q \leq q_1 \leq 2Q, \, 1 \leq a_1 \leq q_1^2 \,\,\textrm{such that }  q_1 \textrm{ prime and the following holds } aq_1^2-a_1q^2 \leq Q.$$
  
  % Zhao has a bound of size $Q^{1/2+\epsilon}$. 

  \section*{Acknowledgements}

The author gratefully acknowledges comments from Karin Halupczok, specifically Remark $1.4$. The author is supported by the Austrian Science Fund (FWF) project Y-901 headed by Christoph Aistleitner.

\end{document}